\documentclass{amsart}

\usepackage{amssymb,amsmath,amscd,graphicx,amsthm}
\newtheorem{theorem}{Theorem}
\newtheorem{lemma}[theorem]{Lemma}

\newtheorem{corollary}[theorem]{Corollary}
\newtheorem{conjecture}[theorem]{Conjecture}

\theoremstyle{definition}

\theoremstyle{remark}
\newtheorem{remark}[theorem]{Remark}

\numberwithin{equation}{section}
\numberwithin{theorem}{section}

\newcommand {\Z} {\mathcal Z}

\newcommand {\CC} {\mathcal C}

\newcommand{\PP}{\mathbb P}

\def\qed{$\;\;\square$}

\def\A{{\mathcal A}}
\def\AA{{\mathbb A}}

\def\C{{\mathbb C}}
\def\CC{{\mathcal C}}

\def\FF{{\mathcal F}}
\def\FFF{{\mathcal F}}
\def\G{{\mathcal G}}
\def\H{\mathcal H}

\def\M{{\mathcal M}}

\def\O{{\mathcal O}}
\def\P{{\mathcal P}}
\def\PP{{\mathcal P}}

\def\TE{{\mathcal T}}

\def\UU{\overline{\A}}

\def\X{{\mathcal X}}
\def\Y{{\mathcal Y}}
\def\Z{{\mathbb Z}}

\def\d{\mathbf d}

\def\g{\mathfrak g}

\def\h{\mathfrak h}

\def\one{\mathbf 1}

\def\phhi{{\varphi}}
\def\pssi{{\psi}}

\def\thetta{{\theta}}

\def\wB{{\widetilde{B}}}
\def\wM{{\widetilde{M}}}

\def\hM{{\widehat{M}}}

\def\wx{{\widetilde{\bf x}}}
\def\x{{\bf x}}

\def\End{\operatorname{End}}

\def\Jac{{\operatorname {Jac}}}

\def\Mat{\operatorname{Mat}}

\def\Poi{{\{\cdot,\cdot\}}}

\def\ord{\operatorname{ord}}

\def\rank{\operatorname{rank}}

\def\:{{:\ }}

\begin{document}

%%%%%%%%%%%%%%%%%%%%%%%%%%%%%%

%% For titles, only capitalize the first letter
%% \title{Almost sharp fronts for the surface quasi-geostrophic equation}

\title{Cremmer--Gervais cluster structure on $SL_n$}

%% Enter authors via the \author command.  
%% Use \affil to define affiliations.
%% (Leave no spaces between author name and \affil command)

%% Note that the \thanks{} command has been disabled in favor of
%% a generic, reserved space for PNAS publication footnotes.

%% \author{<author name>
%% \affil{<number>}{<Institution>}} One number for each institution.
%% The same number should be used for authors that
%% are affiliated with the same institution, after the first time
%% only the number is needed, ie, \affil{number}{text}, \affil{number}{}
%% Then, before last author ...
%% \and
%% \author{<author name>
%% \affil{<number>}{}}

%% For example, assuming Garcia and Sonnery are both affiliated with
%% Universidad de Murcia:
%% \author{Roberta Graff\affil{1}{University of Cambridge, Cambridge,
%% United Kingdom},
%% Javier de Ruiz Garcia\affil{2}{Universidad de Murcia, Bioquimica y Biologia
%% Molecular, Murcia, Spain}, \and Franklin Sonnery\affil{2}{}}
\author{M. Gekhtman}
\address{Department of Mathematics, University of Notre Dame, Notre Dame,
IN 46556}
\email{mgekhtma@nd.edu}

\author{M. Shapiro}
\address{Department of Mathematics, Michigan State University, East Lansing,
MI 48823}
\email{mshapiro@math.msu.edu}

\author{A. Vainshtein}
\address{Department of Mathematics \& Department of Computer Science, University of Haifa, Haifa,
Mount Carmel 31905, Israel}
\email{alek@cs.haifa.ac.il}

\begin{abstract} 
We study natural cluster structures in the rings of regular functions on simple complex Lie groups and Poisson--Lie
structures compatible with these cluster structures.
According to our main conjecture, each class in the Belavin--Drinfeld classification of Poisson--Lie structures on $\G$ corresponds to a cluster structure in $\O(\G)$. 
We have shown before that this conjecture holds for any $\G$ in the case of the standard Poisson--Lie structure   and for all Belavin-Drinfeld classes in $SL_n$, $n<5$.
In this paper we establish it for the Cremmer--Gervais
Poisson--Lie structure on $SL_n$, which is the least similar to the standard one. 
Besides, we prove that on $SL_3$ the cluster algebra and
the upper cluster algebra corresponding to the Cremmer--Gervais cluster structure do not coincide, unlike the case of the standard cluster structure. Finally, we show that the positive locus with respect to the 
 Cremmer--Gervais cluster structure is contained in the set of totally positive matrices.
 \end{abstract}


%% The \maketitle command is necessary to build the title page.
\maketitle

%%%%%%%%%%%%%%%%%%%%%%%%%%%%%%%%%%%%%%%%%%%%%%%%%%%%%%%%%%%%%%%%

\section{Introduction}

In \cite{GSVMMJ} we initiated a systematic study of
multiple cluster structures in the rings of regular functions
on simple Lie groups following an approach developed and implemented 
in \cite{GSV1, GSV5, GSVb} 
for constructing cluster structures 
on algebraic Poisson varieties. It is based on the notion of a Poisson bracket compatible with a cluster structure. The key point is that if an algebraic 
Poisson variety $\left ( \mathcal{M}, \Poi\right )$ possesses a  coordinate chart
that consists of regular functions whose logarithms have pairwise constant Poisson brackets, 
then one can try to use this chart to define a cluster structure $\CC_\M$ compatible with $\Poi$. Algebraic 
structures corresponding to $\CC_\M$ (the cluster algebra and the upper cluster algebra)
are closely related 
to the ring $\O(\M)$ of regular functions on $\mathcal{M}$. 
In fact, under certain rather mild conditions, $\O(\M)$ can be obtained by tensoring one of these
algebras with $\C$.

This construction was applied in \cite[Ch.~4.3]{GSVb} to double Bruhat cells in semisimple Lie groups
equipped with (the restriction of) the {\em standard\/} Poisson--Lie structure. It was shown that
the resulting cluster structure coincides with the one built in \cite{CAIII}. Recall that it was proved in
\cite{CAIII} that the corresponding upper cluster algebra is isomorphic to the ring of regular functions on 
the double Bruhat cell. Since the open double Bruhat cell is dense in the corresponding Lie group, 
one can equip  the ring of regular 
functions on the Lie group with the same cluster structure. 
The standard Poisson--Lie structure is a particular case of Poisson--Lie structures corresponding to quasi-triangular
Lie bialgebras. Such structures are associated with solutions to the classical Yang--Baxter equation (CYBE).
Their complete classification was obtained by Belavin and Drinfeld in \cite{BD}. In \cite{GSVMMJ} we conjectured that any such solution
 gives rise to a compatible cluster structure on the Lie group
 and provided several examples supporting this conjecture by showing that it holds true for the class of the standard Poisson--Lie
  structure in any simple complex Lie group, and for 
 the whole Belavin--Drinfeld classification in $SL_n$ for $n=2,3,4$.  We call the cluster structures
 associated with the non-trivial Belavin--Drinfeld data {\it exotic}.

In this paper, we outline the proof of the conjecture of \cite{GSVMMJ} in the case of the Cremmer--Gervais Poisson structure on $SL_n$. We chose to consider this case because, in a sense, the Poisson structure in question differs the most from the standard: the discrete data featured in the Belavin--Drinfeld classification is trivial in the standard case and has the ``maximal size'' in the Cremmer--Gervais case. Our result allows to equip $SL_n$, $GL_n$ and the affine space $\Mat_n$ of $n\times n$ matrices with a new cluster structure, $\CC_{CG}$. 

In the first Section  below, we collect the necessary information on cluster algebras, compatible Poisson 
brackets, and the toric action. In the next Section, we
%and Poisson--Lie groups and 
formulate the main conjecture from \cite{GSVMMJ}, present the definition of the Cremmer--Gervais Poisson bracket, and formulate our main result. We introduce the cluster structure $\CC_{CG}$ and outline the proof of the main theorem by breaking it into a series of intermediate results about $\CC_{CG}$.
%and showing how the final part of the proof --- the fact that the ring of regular functions is contained in the 
%upper cluster algebra associated with $\CC_{CG}$ --- follows from them. 
In the following Section, we discuss the relation between cluster algebras and upper cluster algebras on $SL_n$. In the standard case these two algebras coincide. We show that for the Cremmer--Gervais cluster structure on $SL_3$ this is not the case. The next Section treats  positivity for the exotic cluster structure $\CC_{CG}$.
Finally, in the last Section we formulate several directions for future research.
%there exists a gap between the two.
 
\section{Cluster structures and compatible Poisson brackets}
\label{SecPrel}

We start with the basics on cluster algebras of geometric type. The definition that we present
below is not the most general one, see, e.g.,
\cite{FZ2, CAIII} for a detailed exposition. In what follows, we will use a notation $[i,j]$ for an interval
$\{i, i+1, \ldots , j\}$ in $\mathbb{N}$, and we will denote $[1, n]$ by $[n]$.
 
The {\em coefficient group\/} $\PP$ is a free multiplicative abelian
group of finite rank $m$ with generators $g_1,\dots, g_m$.
An {\em ambient field\/}  is
the field $\FFF$ of rational functions in $n$ independent variables with
coefficients in the field of fractions of the integer group ring
$\Z\PP=\Z[g_1^{\pm1},\dots,g_m^{\pm1}]$ (here we write
$x^{\pm1}$ instead of $x,x^{-1}$).

A {\em seed\/} (of {\em geometric type\/}) in $\FFF$ is a pair
$\Sigma=(\x,\widetilde{B})$,
where $\x=(x_1,\dots,x_n)$ is a transcendence basis of $\FFF$ over the field of
fractions of $\Z\PP$ and $\widetilde{B}$ is an $n\times(n+m)$ integer matrix
whose principal part $B$ is skew-symmetrizable (recall that the principal part of a rectangular matrix  
is its maximal leading square submatrix).  Matrices $B$ and $\wB$ are called the
{\it exchange matrix\/} and the {\it extended exchange matrix}, respectively.
In this paper, we will only deal with the case when the exchange matrix is skew-symmetric.

The $n$-tuple  $\x$ is called a {\em cluster\/}, and its elements
$x_1,\dots,x_n$ are called {\em cluster variables\/}. Denote
$x_{n+i}=g_i$ for $i\in [m]$. We say that
$\widetilde{\x}=(x_1,\dots,x_{n+m})$ is an {\em extended
cluster\/}, and $x_{n+1},\dots,x_{n+m}$ are {\em stable
variables\/}. It is convenient to think of $\FFF$ as
of the field of rational functions in $n+m$ independent variables
with rational coefficients. 

In what follows, we will only deal with the case when the exchange matrix is skew-symmetric. In this
situation the extended exchange  
matrix can be conveniently represented by a {\it quiver\/} $Q=Q(\wB)$. It is a directed graph on the 
vertices $1,\dots,n+m$ corresponding to all variables; the vertices corresponding to stable variables are called
{\it frozen}. Each entry $b_{ij}>0$ of the matrix $\wB$ gives rise to $b_{ij}$ edges going from
the vertex $i$ to the vertex $j$; each such edge is denoted $i\to j$. Clearly, $\wB$ can be restored uniquely from $Q$. 
%and we will eventually write $\Sigma=(\x,Q)$ instead of $\Sigma=(\x,\wB)$.
%Note that $B$ is irreducible if and only if the 
%subquiver of $Q$ induced by non-frozen vertices is connected.

Given a seed as above, the {\em adjacent cluster\/} in direction $k\in [n]$
is defined by
$$
\x_k=(\x\setminus\{x_k\})\cup\{x'_k\},
$$
where the new cluster variable $x'_k$ is given by the {\em exchange relation}
\begin{equation}\label{exchange}
x_kx'_k=\prod_{\substack{1\le i\le n+m\\  b_{ki}>0}}x_i^{b_{ki}}+
       \prod_{\substack{1\le i\le n+m\\  b_{ki}<0}}x_i^{-b_{ki}};
\end{equation}
here, as usual, the product over the empty set is assumed to be
equal to~$1$.

We say that $\wB'$ is
obtained from $\wB$ by a {\em matrix mutation\/} in direction $k$
and
write $\wB'=\mu_k(\wB)$ 
 if
\[
b'_{ij}=\begin{cases}
         -b_{ij}, & \text{if $i=k$ or $j=k$;}\\
                 b_{ij}+\displaystyle\frac{|b_{ik}|b_{kj}+b_{ik}|b_{kj}|}2,
                                                  &\text{otherwise.}
        \end{cases}
\]
It can be easily verified that $\mu_k(\mu_k(\wB))=\wB$.

Given a seed $\Sigma=(\x,\widetilde{B})$, we say that a seed
$\Sigma'=(\x',\widetilde{B}')$ is {\em adjacent\/} to $\Sigma$ (in direction
$k$) if $\x'$ is adjacent to $\x$ in direction $k$ and $\widetilde{B}'=
\mu_k(\widetilde{B})$. Two seeds are {\em mutation equivalent\/} if they can
be connected by a sequence of pairwise adjacent seeds. 
The set of all seeds mutation equivalent to $\Sigma$ is called the {\it cluster structure\/} 
(of geometric type) in $\FFF$ associated with $\Sigma$ and denoted by $\CC(\Sigma)$; in what follows, 
we usually write $\CC(\wB)$, or even just $\CC$ instead. 

Following \cite{FZ2, CAIII}, we associate
with $\CC(\wB)$ two algebras of rank $n$ over the {\it ground ring\/} $\AA$, $\Z\subseteq\AA \subseteq\Z\P$:
the {\em cluster algebra\/} $\A=\A(\CC)=\A(\wB)$, which 
is the $\AA$-subalgebra of $\FF$ generated by all cluster
variables in all seeds in $\CC(\wB)$, and the {\it upper cluster algebra\/}
$\UU=\UU(\CC)=\UU(\wB)$, which is the intersection of the rings of Laurent polynomials over $\AA$ in cluster variables
taken over all seeds in $\CC(\wB)$. The famous {\it Laurent phenomenon\/} \cite{FZ3}
claims the inclusion $\A(\CC)\subseteq\UU(\CC)$. The natural choice of the ground ring for the geometric type
is the polynomial ring in stable variables $\AA=\Z\P_+=\Z[x_{n+1},\dots,x_{n+m}]$; this choice is assumed unless
explicitly stated otherwise. 

Let $V$ be a quasi-affine variety over $\C$, $\C(V)$ be the field of rational functions on $V$, and
$\O(V)$ be the ring of regular functions on $V$. Let $\CC$ be a cluster structure in $\FF$ as above.
Assume that $\{f_1,\dots,f_{n+m}\}$ is a transcendence basis of $\C(V)$. Then the map $\varphi: x_i\mapsto f_i$,
$1\le i\le n+m$, can be extended to a field isomorphism $\varphi: \FF_\C\to \C(V)$,  
where $\FF_\C=\FF\otimes\C$ is obtained from $\FF$ by extension of scalars.
The pair $(\CC,\varphi)$ is called a cluster structure {\it in\/}
$\C(V)$ (or just a cluster structure {\it on\/} $V$), $\{f_1,\dots,f_{n+m}\}$ is called an extended cluster in
 $(\CC,\varphi)$.
Sometimes we omit direct indication of $\varphi$ and say that $\CC$ is a cluster structure on $V$. 
A cluster structure $(\CC,\varphi)$ is called {\it regular\/}
if $\varphi(x)$ is a regular function for any cluster variable $x$. 
The two algebras defined above have their counterparts in $\FF_\C$ obtained by extension of scalars; they are
denoted $\A_\C$ and $\UU_\C$.
If, moreover, the field isomorphism $\varphi$ can be restricted to an isomorphism of 
$\A_\C$ (or $\UU_\C$) and $\O(V)$, we say that 
$\A_\C$ (or $\UU_\C$) is {\it naturally isomorphic\/} to $\O(V)$.

Let $\Poi$ be a Poisson bracket on the ambient field $\FFF$, and $\CC$ be a cluster structure in $\FFF$. 
We say that the bracket and the cluster structure are {\em compatible\/} if, for any extended
cluster $\widetilde{\x}=(x_1,\dots,x_{n+m})$,  one has
\begin{equation}\label{cpt}
\{x_i,x_j\}=\omega_{ij} x_ix_j,
\end{equation}
where $\omega_{ij}\in\Z$ are
constants for all $i,j\in[n+m]$. The matrix
$\Omega^{\widetilde \x}=(\omega_{ij})$ is called the {\it coefficient matrix\/}
of $\Poi$ (in the basis $\widetilde \x$); clearly, $\Omega^{\widetilde \x}$ is
skew-symmetric. The notion of compatibility  extend to Poisson brackets on $\FF_\C$ without any changes.

A complete characterization of Poisson brackets compatible with a given cluster structure $\CC=\CC(\wB)$ in the case $\rank\wB=n$ is given in \cite{GSV1}, see also \cite[Ch.~4]{GSVb}. 
A different description of compatible Poisson brackets on $\FF_\C$ is based on the notion of a toric action.
 Fix an arbitrary extended cluster
$\wx=(x_1,\dots,x_{n+m})$ and define a {\it local toric action\/} of rank $r$ as the map 
$\TE^W_{\d}:\FF_\C\to
\FF_\C$ given on the generators of $\FF_\C=\C(x_1,\dots,x_{n+m})$ by the formula 
\begin{equation}
\TE^W_{\d}(\wx)=\left ( x_i \prod_{\alpha=1}^r d_\alpha^{w_{i\alpha}}\right )_{i=1}^{n+m},\qquad
\d=(d_1,\dots,d_r)\in (\C^*)^r,
\label{toricact}
\end{equation}
where $W=(w_{i\alpha})$ is an integer $(n+m)\times r$ {\it weight matrix\/} of full rank, and extended naturally to the whole $\FF_\C$. 

Let $\wx'$ be another extended cluster, then the corresponding local toric action defined by the weight matrix $W'$
is {\it compatible\/} with the local toric action \eqref{toricact} if it commutes with the cluster transformation that takes $\wx$ to $\wx'$.
 If local toric actions at all clusters are compatible, they define a {\it global toric action\/} $\TE_{\d}$ on $\FF_\C$ called the extension of the local toric action \eqref{toricact}. 
 If $\rank\wB=n$, then the maximal possible rank of a global toric action equals $m$. Any global toric action can be obtained from a toric action of
the maximal rank by setting some of $d_i$'s equal to~$1$.

\section{Main results}
\label{SecMC}

Let $\G$ be a Lie group equipped with a Poisson bracket $\Poi$.
$\G$ is called a {\em Poisson--Lie group\/}
if the multiplication map
$$
\G\times \G \ni (x,y) \mapsto x y \in \G
$$
is Poisson. The tangent Lie algebra $\g$ of a Poisson--Lie group $\G$ has a natural structure of a {\em Lie bialgebra}. We are interested in the case when
$\G$ is a simple complex Lie group and its tangent Lie bialgebra is {\em factorizable}. 

A factorizable Lie bialgebra structure on  a complex simple Lie algebra can be described in terms of  a classical R-matrix.  This is an element
$r\in \g \otimes \g$ that satisfies the {\em classical Yang-Baxter equation\/} and an additional condition 
that $r + r^{21}$  defines an invariant nondegenerate inner product on $\g$. 
(Here $r^{21}$ is obtained from $r$ by switching factors in tensor products.)  Classical R-matrices were 
classified, up to an automorphism,  by Belavin and Drinfeld in \cite{BD}. 
Let $\h$ be a Cartan subalgebra of $\g$,
$\Phi$ be the root system associated with $\g$, $\Phi^+$ be the set of positive roots, and $\Delta\subset \Phi^+$ be the 
set of positive simple roots. The Killing form on $\g$ is denoted by $\langle\ , \ \rangle$.
A {\em Belavin--Drinfeld triple\/} $T=(\Gamma_1,\Gamma_2, \gamma)$
consists of two subsets $\Gamma_1,\Gamma_2$ of $\Delta$ and an isometry $\gamma:\Gamma_1\to\Gamma_2$ nilpotent in 
the following sense: for every $\alpha \in \Gamma_1$ there exists $m\in\mathbb{N}$ such that $\gamma^j(\alpha)\in 
\Gamma_1$ for $j=0,\ldots,m-1$, but $\gamma^m(\alpha)\notin \Gamma_1$. The isometry $\gamma$ extends in a natural way to a map between root systems $\Phi_1, \Phi_2$ generated by $\Gamma_1, \Gamma_2$. This allows one to define a partial ordering on $\Phi$: $\alpha \prec_T \beta$ if $\beta=\gamma^j(\alpha)$ for some $j\in \mathbb{N}$, and to set $\h_T=\{ h\in\h \:  \alpha(h)=\beta(h)\ \text{if}\ \alpha\prec_T\beta\}$.

To each $T$ there corresponds a set 
$\mathcal R_T$ of classical R-matrices that we  call the  {\em Belavin-Drinfeld class\/} corresponding to $T$. 
Two R-matrices in  the same Belavin-Drinfeld class $\mathcal R_T$ differ by an element from $\h\otimes \h$ 
satisfying a linear relation specified by $T$. We denote by $\Poi_r$ the Poisson--Lie bracket associated with 
$r\in \mathcal R_T$. Given a Belavin--Drinfeld triple $T$ for $\G$,
define the torus $\H_T=\exp \h_T\subset\G$. Note that the dimension of the torus equals $k_T=|\Delta\setminus\Gamma_1|$.

In \cite{GSVMMJ} we conjectured that there exists a classification of regular cluster structures on $\G$ that is 
completely parallel to the Belavin--Drinfeld classification.

\begin{conjecture}
\label{ulti}
Let $\G$ be a simple complex Lie group.
For any Belavin-Drinfeld triple $T=(\Gamma_1,\Gamma_2,\gamma)$ there exists a cluster structure
$(\CC_T,\varphi_T)$ on $\G$ such that

{\rm (i)} %the number of variables in each extended cluster of $\CC_T$ is $\dim \G$, 
the number of stable variables is $2k_T$, and the corresponding extended exchange matrix has a full rank; 
%$2\dim \H_\tau$;

{\rm (ii)} $(\CC_T,\varphi_T)$ is regular, and the corresponding upper cluster algebra $\UU_\C(\CC_T)$ 
is naturally isomorphic to $\O(\G)$;

{\rm (iii)} the global toric action of $(\mathbb{C}^*)^{2k_T}$ on $\C(\G)$ is generated by the action
of $\H_T\times \H_T$ on $\G$ given by $(H_1, H_2)(X) = H_1 X H_2$;

 {\rm (iv)} for any $r\in \mathcal R_T$, $\Poi_r$ is compatible with $\CC_T$;

{\rm (v)} a Poisson--Lie bracket on $\G$ is compatible with $\CC_T$ only if it is a scalar multiple of
 $\Poi_r$ for some $r\in \mathcal R_T$.
\end{conjecture}

The Belavin-Drinfeld data (triple, class) is said to be {\it trivial\/} if $\Gamma_1=\Gamma_2=\varnothing$.
In this case,   $\H_T=\H$ is the Cartan
subgroup in $\G$. 
The resulting Poisson bracket is called  {\em the standard Poisson--Lie structure} on $\G$.
The Conjecture \ref{ulti} in this case was verified in \cite{GSVMMJ}.

In this paper we consider the case $\G=SL_n$ and the Belavin--Drinfeld data that is "the farthest" from the trivial data, namely, $\Gamma_1=\{\alpha_2, \dots, \alpha_{n-1}\}, \ 
 \Gamma_2=\{\alpha_1, \dots, \alpha_{n-2}\}$ and $\gamma(\alpha_i) = \alpha_{i-1}$
for $i=2, \ldots,n-1$.  The resulting Poisson--Lie bracket on $SL_n$ is called the
{\em Cremmer--Gervais bracket}. The main result of this paper is

\begin{theorem}\label{CGtrue}
Conjecture {\rm \ref{ulti}} is valid for the Cremmer--Gervais Poisson--Lie structure.
\end{theorem}

Theorem~\ref{CGtrue} is proved by producing a cluster structure $\CC_{CG}=\CC_{CG}(n)$ that possesses all the
needed properties. In fact, we will construct a cluster structure in the space $\Mat_n$ of 
$n\times n$ matrices compatible with a natural extension of the Cremmer--Gervais Poisson bracket and derive the 
required properties of $\CC_{CG}$ from similar features of the latter cluster structure. 
%Note that in the ``intermediate'' case of $GL_n$ we also obtain a regular cluster structure compatible
%with the extension of the Cremmer--Gervais Poisson bracket, however, in this case the ring of regular functions on $GL_n$ is isomorphic to the localization of the upper cluster algebra with respect to the function $\det X$.
In what follows we use the same notation $\CC_{CG}$ for both cluster structures and indicate explicitely which
one is meant when needed.

To describe the initial cluster for  $\CC_{CG}$, we need to introduce some notation.
For a matrix $A$, we denote by $A_{i_1\ldots i_l}^{j_1\ldots j_m}$ its submatrix formed
by rows $i_1,\ldots, i_l$ and columns $j_1,\ldots, j_m$. If all rows (respectively, columns) of $A$ are 
involved, we will omit the lower (respectively, upper) list of indices.
If $X$, $Y$ are two $n\times n$ matrices, denote by  $\X$ and $\Y$  $(n-1)\times (n+1)$ matrices 
$$
\X = \left [ X_{[2,n]}\  0\right ], \quad \Y = \left [0\   Y_{[1,n-1]} \right ].
$$
Put $k=\lfloor \frac{n+1}{2}\rfloor$ and $N=k (n-1)$. Define a $k (n-1) \times (k+1) (n+1)$ matrix
\begin{equation}
U(X, Y) = \left [
\begin{array}{ccccc}
\Y & \X & 0 & \cdots & 0\\
0 & \Y & \X  & 0 & \cdots\\
0 & \ddots& \ddots &\ddots & 0\\
0 & \cdots & 0 & \Y & \X
\end{array}
\right ].
\label{uho}
\end{equation}

Define three families of functions %in $X,Y\in SL_n$ 
via
\begin{equation}\label{inclust}
 \begin{aligned}
\thetta_i(X)&=\det X_{[n-i+1,n]}^{[n-i+1,n]}, \; i\in [n-1];\\
\phhi_p(X,Y)&=\det U(X,Y)_{[N-p+1, N]}^{[k (n+1) - p +1, k (n+1)]}, \; p\in [N];\\
\pssi_q(X,Y)&=\det U(X,Y)_{[N-q+1, N]}^{[k (n+1) - q +2, k (n+1)+1]}, \; q\in [M].
\end{aligned}
\end{equation}
In the last family, $M=N$ if $n$ is even and $M=N-n+1$ if $n$ is odd.

\begin{theorem}
\label{logcan}
The functions
$\theta_i(X)$, $\phi_p(X,X)$, $\psi_q(X,X)$ form a log-canonical family with respect to the Cremmer--Gervais
bracket. 
\end{theorem}

Consequently, we choose family~\eqref{inclust} resricted to $X=Y$ 
as an initial (extended) cluster for  $\CC_{CG}(n)$. 
Furthermore, functions $\phhi_N$ and $\pssi_M$ are the only stable variables for $\CC_{CG}(n)$.

Let us explain the motivation behind this choice and the intuition it provides for constructing the initial cluster for $\CC_{CG}(n)$ presented above. 
To this end, we first need to recall the notion of the {\em Drinfeld double\/} $D(\G)$ of a factorizable Poisson--Lie group $(\G, \Poi_r)$. The double 
$D(\G)$ is endowed with a Poisson--Lie structure associated with the {\em Manin triple\/} $(D(\g), \g_d, \g_r)$. 
Here $D(\g)=\g \oplus \g = \g_d \dot + \g_r$ is equipped with the invariant nondegenerate inner product 
$\langle\langle (\xi,\eta), (\xi',\eta')\rangle\rangle = \langle \xi, \xi'\rangle - \langle \eta, \eta'\rangle$, an isotropic subalgebra $\g_d$ is the image
of $\g$ in $D(\g)$ under the diagonal embedding, and $\g_r$ is an isotropic subalgebra of $D(\g)$ 
given by $\g_r=\{(R_+(\xi),R_-(\xi))\: \xi \in \g\}$, where $R_\pm\in \End\g$ are defined by $\langle R_+ \eta, \zeta \rangle = - \langle R_- \zeta, \eta \rangle =\langle r, \eta\otimes\zeta \rangle$. The embedding 
$\g\hookrightarrow D(\g)$ whose image is $\g_r$ is denoted by $\sigma$.
The group $(\G,\Poi_r)$  becomes a Poisson--Lie subgroup of $D(\G)$ under the diagonal embedding.
Another Poisson--Lie subgroup of $D(\G)$
is the group $\G_r$  whose Lie algebra is $\g_r$. 
Orbits of the two-sided action of $\G_r \times \G_r$ on $D(\G)$ play an important role in the description of symplectic leaves of  $D(\G)$ and $\G$ (\cite{r-sts, Ya}). This action 
is relevant to our tackling on Conjecture \ref{ulti} due to the following observation. Let $f$ be a {\em semi-invariant\/} of the two-sided action of $\G_r \times \G_r$, that is, for any $g_1, g_2 \in \G_r$ 
and $(X,Y) \in D(G)$, $f(g_1 (X,Y) g_2)= \xi_l(g_1) f(X,Y) \xi_r(g_2)$, where $\xi_l, \xi_r$ are two characters of $\G_r$ trivial on elements of the form $\exp(\sigma(e_\alpha))$ for any root vector $e_\alpha$.
Then the Hamiltonian flow generated by $f(X,X)$ on $\G$ with respect to $\Poi_r$ is given by
$X(t) = \exp(t h_l) X(0) \exp( t h_r)$, where $h_l, h_r$ are two elements of the Cartan subalgebra.

This observation makes  regular functions on $\G$ obtained by restriction to the diagonal subgroup of semi-invariants described above natural candidates for stable variables in the cluster structure we are trying to 
construct. Indeed,
if $y$ is a  stable variable in a cluster algebra admitting a compatible Poisson bracket $\Poi$, then, for any cluster variable
$x$, $\{\log y,x\}= c_x x$ for some $c_x\in {\mathbb Z}$, and thus  the Hamiltonian flow generated
by $\log y$ gives rise to a global toric action $\mathcal T_t(x)= x \exp(c_x t)$. On the other hand, part (iii) of Conjecture~\ref{ulti} predicts that the global toric action in $\CC_{CG}(n)$ is induced by 
the right and left action on $\G$ by elements of the subgroup $\mathcal H_T$ of the Cartan group specified by the Belavin--Drinfeld data. 

The stable variables for $\CC_{CG}(n)$ were constructed using the strategy indicated above: the action $(X,Y) 
\mapsto g_1 (X,Y) g_2$ translates into a transformation $U(X,Y)\mapsto M_l(g_1) U(X,Y) M_r(g_2)$, where  $M_l(g_1), M_r(g_2)$ are certain invertible block diagonal matrices. This allows to identify $\phhi_N(X,Y)$ and $\pssi_M(X,Y)$ as semi-invariants, and also hints at a possibility 
of constructing an initial cluster consisting of minors of $U(X,Y)$. It is this intuition that eventually led us to formulation of Theorem \ref{logcan}. Moreover, the proof of the theorem relies on computation of Poisson brackets in the double.

Next, we need to describe the quiver  that corresponds to the initial cluster above. 
In fact, it will be convenient to do this in the $\Mat_n$ rather than $SL_n$ situation, in
which case the initial cluster is augmented by the  addition of one more stable
variable, 
$\thetta_n(X)=\det X$. The Cremmer--Gervais Poisson structure is extended to $\Mat_n$
by requiring that $\thetta_n(X)$ is a Casimir function.
We denote the quiver corresponding to the augmented initial cluster by $Q_{CG}(n)$.
Its vertices are $n^2$ nodes of the $n\times n$ rectangular grid indexed by pairs
$(i,j)$, $i, j \in [n]$, with $i$ increasing top to bottom and $j$ increasing left
to right. Before
describing  edges of the quiver, let us explain the correspondence between the
cluster variables and the vertices of $Q_{CG}(n)$.

For any cluster or stable variable $f$ in the augmented initial cluster, consider the upper
left matrix
entry of the submatrix of $U(X,X)$ (or $X$) associated with this variable. This
matrix entry is $x_{ij}$ for some $i, j \in [n]$. 
Thus we define a correspondence  $\rho \: f \leftrightarrow (i,j)\in
[n]\times [n]$. We have proved that
$\rho$ is a one-to-one correspondence between the augmented initial cluster and
$[n]\times [n]$.

Let us assign each cluster variable $f$ to the vertex indexed by $(i,j) = \rho(f)$.
In particular, the stable variable $\thetta_n$ is assigned to vertex $(1,1)$, 
 the stable variable $\phhi_N$ is assigned to vertex $(2,1)$ if $n$ is odd and $(1,n)$
 if $n$ is even, and  the stable variable $\pssi_M$ is assigned to vertex $(1,n)$ if
$n$ is odd and $(2,1)$
 if $n$ is even.

Now, let us describe arrows of $Q_{CG}(n)$. There are horizontal arrows
$(i,j+1) \to (i,j)$ for  all $i \in [n], j\in [n-1]$ except $(i,j)=(1, n-1)$;
vertical arrows
$(i+1,j) \to (i,j)$ for  all $i \in [n-1], j\in [n]$ except $(i,j)=(1,1)$; 
diagonal arrows $(i,j) \to (i+1,j+1)$ for  all $i, j\in [n-1]$.
In addition, there are arrows between vertices of the first and the last rows:
$(n,j) \to (1,j)$ and $(1,j) \to (n,j+1)$ for $j\in [2, n-1]$;  and 
 arrows between vertices of the first and the last columns:
$(i,n) \to (i+2,1)$ and $(i+2,1) \to (i+1, n)$ for $i\in [1, n-2]$. This
concludes the description of $Q_{CG}(n)$. The quiver  $Q'_{CG}(n)$ that correspond to the
$SL_n$ case is obtained
from $Q_{CG}(n)$ by deleting the vertex $(1,1)$ and erasing all arrows incident to this
vertex. Quiver $Q_{CG}(5)$ is shown on Fig.~1. 

\begin{figure}[ht]
\begin{center}
\includegraphics[height=6.4cm]{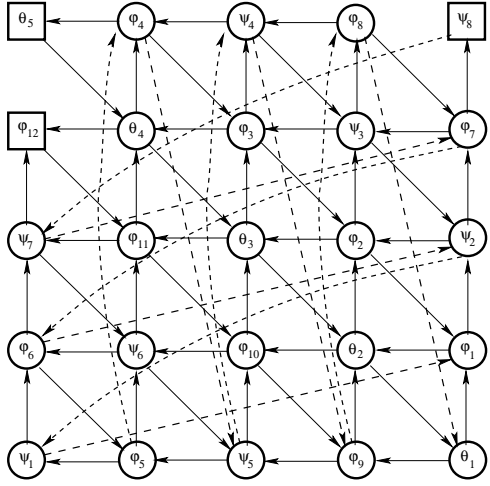}
\caption{Quiver$Q_{CG}(5)$}
\label{fig:Qcg45}
\end{center}
\end{figure}

\begin{theorem}
\label{quiver} {\rm (i)}
The quivers $Q_{CG}(n)$  and $Q'_{CG}(n)$ define cluster structures compatible with the
Cremmer--Gervais Poisson structure on $\Mat_n$ and $SL_n$ respectively.

{\rm (ii)} The corresponding extended exchange matrices are of full rank.
\end{theorem}

Note that $k_T=1$ in the Cremmer--Gervais case, and hence the corresponding
Belavin--Drinfeld class contains a unique R-matrix. Therefore, Theorem~\ref{quiver} establishes
parts (i) and (iv) of Conjecture~\ref{ulti}  in the Cremmer--Gervais case.

Another property of the cluster structure $\CC_{CG}$ is given by the following theorem.

\begin{theorem}
\label{regular}
The  cluster structure $\CC_{CG}$ is regular.
\end{theorem}

The proof of Theorem~\ref{regular} relies on  Dodgson-type identities applied to
submatrices of $U(X,Y)$ while taking into account its shift-invariance properties. 
As a corollary, we get parts (iii) and (v) of  Conjecture~\ref{ulti}.

By Theorem~\ref{regular}, the upper cluster algebra $\UU_{CG}=\UU_{\C}(\CC_{CG})$ is a
subalgebra in the ring of regular
 functions. To complete the proof of Theorem~\ref{CGtrue}, it remains to
establish the opposite inclusion, which will settle part (ii) of Conjecture~\ref{ulti}.
 The proof relies on induction on $n$. Its main ingredient is a construction of
two distinguished sequences of cluster
 transformations. The first sequence, $\mathcal S$, followed by freezing some of the
cluster variables and localization
 at a single cluster variable $\phhi_{n-1}(X)$, leads to a map $\zeta \: \Mat_n \setminus
 \{X\:\phhi_{n-1}(X)=0\}\to
\Mat_{n-1}$ that ``respects" the Cremmer--Gervais cluster
 structure. The map $\zeta$ is needed to perform an induction step. However, because
of the localization mentioned above, we also need a second sequence, $\mathcal T$, 
of transformations that can be viewed as a cluster-algebraic realization of the
anti-Poisson involution
 $X \mapsto W_0 X W_0$ on $\Mat_n$ equipped with the Cremmer--Gervais Poisson bracket. This allows one to apply $\zeta$ to
  $W_0 X W_0$ as well and then invoke certain general properties of cluster algebras
to fully utilize the induction assumption.

Denote by $\hat Q_{CG}(n)$ the quiver obtained by adding to $Q_{CG}(n)$ two
additional arrows: $(1,1) \to (n,2)$ and $(n,1)\to (1,1)$.

 \begin{theorem}
 \label{transform1}
  There exists a sequence $\mathcal S$ of cluster transformations in $\CC_{CG}(n)$
such that $\mathcal S(Q_{CG}(n))$ contains
  a subquiver isomorphic to $\hat Q_{CG}(n-1)$ and cluster variables indexed by
vertices of this subquiver satisfy
  $$\mathcal S(f)_{ij}(X)=
  \phhi_{n-1}(X)^{\varepsilon_{ij}}
  f_{ij}(\zeta(X)),\quad i,j\in [n-1],
  $$
  where $\varepsilon_{ij}=1$ if $(i,j)$ is a $\pssi$-vertex in $Q_{CG}(n-1)$
and $\varepsilon_{ij}=0$ otherwise.
 \end{theorem}
 
For any function $g$ on $\Mat_n$ define $g^{w_0}(X)=g(W_0XW_0)$.
 Besides, for any quiver $Q$ denote by $Q^{opp}$ the quiver obtained from $Q$ by reversing all
arrows. 

  \begin{theorem}
 \label{transform2}
 There exists a sequence $\mathcal T$ of cluster transformations in $\CC_{CG}$ such
that $\mathcal T(Q_{CG}(n))$ is isomorphic
 to $Q_{CG}^{opp}(n)$ and cluster variables indexed by vertices of
 $\mathcal T(Q_{CG}(n))$ satisfy
 $$\mathcal T(f)_{ij}(X) = f^{w_0}_{ij}(X)$$.
 \end{theorem}
 
Each cluster transformation in the sequences $\mathcal S$ and $\mathcal T$
 can be written as a Dodgson-type identity for minors of an augmentation of the matrix $U(X,X)$. 

 To complete the proof of Theorem \ref{CGtrue}, it remains to show that functions $x_{ij} \in \O (\Mat_n)$ belong 
 to the upper cluster algebra $\UU_{CG}(n)$, which is done by induction on $n$ that relies on
 Theorems~\ref{transform1} and~\ref{transform2}. The details of the proofs can be found in~\cite{GSV6}.

%\section{The gap between the cluster algebra and the upper cluster algebra}
\section{Relation between $\A_{CG}$ and $\UU_{CG}$}

Recall that the standard cluster algebra on $SL_n$, that is, the one that corresponds to the standard Poisson--Lie structure, coincides with the standard upper cluster algebra. In this section we will prove that this is not the case for the Cremmer--Gervais cluster structure. More exactly, we prove the following theorem.

\begin{theorem}\label{gap}
The cluster algebra $\A_{CG}(3)$ is a proper subalgebra  of the upper cluster algebra $\UU_{CG}(3)$.
\end{theorem}

{\bf Proof:} It will be convenient to consider regular functions on $\Mat_3$,
instead of $SL_3$. 
Recall that the Cremmer--Gervais cluster structure $\CC_{CG}$ in $\Mat_3$ has six cluster variables and 
three stable variables $s_1, s_2, s_3$, where
\begin{equation*}
\begin{gathered}
s_1=\pssi_2(X)=
{\begin{vmatrix}
x_{13} & x_{21} \\
x_{23} & x_{31}
\end{vmatrix},} \\
s_2=\phhi_4(X)=
{\begin{vmatrix}
x_{21} & x_{22} & x_{23} & 0 \\
x_{31} & x_{32} & x_{33} & 0 \\
0 & x_{11} & x_{12} & x_{13} \\
0 & x_{21} & x_{22} & x_{23}
\end{vmatrix}}, \\
 s_3=\thetta_3(X)=|X|.
\end{gathered}
\end{equation*}
Matrix entries $x_{11}$, $x_{12}$, $x_{13}$, $x_{21}$, $x_{22}$, $x_{23}$ together with stable variables $s_1$, $s_2$, 
$s_3$ form a coordinate system on an open subset in $\Mat_3$.
It will be important for the future to calculate the Poisson bracket of the matrix entry $x_{12}$ with
all the other coordinates:
\begin{equation}\label{coorbra}
\begin{aligned}
\{x_{12},x_{11}\} &= -\frac{2}{3} x_{11}x_{12}, \quad
\{x_{12},x_{13}\} =  \frac{2}{3} x_{12}x_{13}, \\
\{x_{12},x_{21}\}& =  \frac{2}{3} x_{12}x_{21}, \quad
\{x_{12},x_{23}\} =   2 x_{13}x_{22}, \\
\{x_{12},x_{22}\} & =  \frac{4}{3} x_{12}x_{22} +2x_{13}x_{21}-2x_{11}x_{23}, \\  
\{x_{12},s_1\}& = \frac{2}{3} x_{12}s_1, \quad  
\{x_{12},s_2\} =  \frac{4}{3} x_{12}s_2, \quad 
\{x_{21},s_3\} =  0.
\end{aligned}
\end{equation}

Consider the polynomial
$$
p= 
\begin{vmatrix}
x_{11} & x_{12} & x_{13} & 0 \\
x_{21} & x_{22} & x_{23} & 0 \\
0 & x_{11} & x_{12} & x_{13} \\
0 & x_{21} & x_{22} & x_{23}
\end{vmatrix}.
$$ 
Note that it can be written as $p=p_0+x_{12}p_1+x_{12}^2p_2$ with
\begin{equation*}
\begin{gathered}
p_0=(x_{11}x_{13}) x_{22}^2+(x_{11} x_{23}-x_{13}x_{21})^2,\\ 
p_1=x_{22}\left(x_{13}x_{21}+x_{11}x_{23}\right), \qquad
p_2=x_{21}x_{23}.
\end{gathered}
\end{equation*}

In what follows we consider the following rings: 
\begin{equation*}
\begin{gathered}
M=\C[x_{11}, x_{13}, x_{21}, s_1, s_2, s_3]_{x_{11}x_{13}x_{21}},\\
 \hM=M[x_{22}, x_{23}]_{p_0}, \qquad \wM=\hM[x_{12}]_p,
\end{gathered}
\end{equation*}
where the subscript stands for the localization. 
Expressing the last row of matrix $X$ through the first two rows and 
$s_1, s_2, s_3$ we observe that 
$\UU_{CG}(3)=\C[\Mat_3]\subset \wM$.

Define two differential operators 
%in coordinates $S_1,S_2,S_3,x_{11},x_{12},x_{13}, x_{21},x_{22},x_{23}$:
\[\begin{aligned}
D_1&=\frac{2}{3}x_{12}\left(-\frac{\partial}{\partial x_{11}}+
\frac{\partial}{\partial x_{13}}+
\frac{\partial}{\partial x_{21}}+
2\frac{\partial}{\partial x_{22}}+
\frac{\partial}{\partial s_1}+
2\frac{\partial}{\partial s_2}\right),\\ 
D_2&=2\left(x_{13}x_{21}-x_{11}x_{23}\right)\frac{\partial}{\partial x_{22}}+
 2 x_{13}x_{22}\frac{\partial}{\partial x_{23}}.
\end{aligned} 
 \]
It follows from \eqref{coorbra} that $\{x_{12},f\}=(D_1+D_2)f$ for any $f\in \wM$.

The ring $\wM$ has a filtration by ideals $\wM\supset x_{12}\wM\supset x_{12}^2 \wM\supset\cdots$.
Clearly,  $\cap_{k=0}^\infty x_{12}^k \wM=0$.
Define for any nonzero $f\in \wM$ the order $\ord f$ as 
the maximal $\nu\ge 0$ such that $f\in x_{12}^\nu \wM$.  Set $\ord 0=\infty$.
Note  that for any $f\in \wM$ the order satisfies inequalities 
\[
\ord D_1(f)\ge \ord f+1,\qquad \ord D_2(f)\ge \ord f.
\]

\begin{lemma}\label{lem:decomp}
Let $f\in \wM$ and $\ord f=0$, then there is a unique decomposition  $f=[f]+x_{12}\tilde f$ 
with $[f]\in \hM$ and $\tilde f \in\wM$.
\end{lemma}

{\bf Proof:}
The existence of the decomposition follows from the identity $1/p=1/p_0-x_{12}(p_1+x_{12}p_2)/(pp_0)$.
 Uniqueness follows from the observation that no element of the subring 
 $\hM$ depends on $x_{12}$. 
\qed

For $f\in\wM$ such that $\ord f>0$ we put $[f]=0$.

\begin{lemma}\label{lem:level}
Let $a, b, c, x, y, u_1,\dots u_m$ be independent variables, and let polynomials
$g=abx^2+(ay-bc)^2$ and $f\in\C[a,b,c;u_1,\dots,u_m;x,y]_{abcg}$ satisfy the following condition: 
for  any choice of values $a=a^0, b=b^0,c=c^0,u_i=u_i^0$ in an open subset of $\C^{m+3}$
and any $\sigma\in\C$ there exists $\tau\in\C$ such that 
\begin{equation*}
\begin{aligned}
\left.\{g\right\vert_{a=a^0, b=b^0, c=c^0, u_i=u_i^0}=\sigma\}&\subset
\left.\{f\right\vert_{a=a^0,b=b^0,c=c^0,u_i=u_i^0}=\tau\}\\
&\subset\C^2.
\end{aligned}
\end{equation*} 
Then $f=F(g)$ where $F$ is a Laurent polynomial 
over $\C[a,b,c;u_1,\dots,u_m]_{abc}$.
\end{lemma}

{\bf Proof:} Multiplying $f$ if needed by a positive power of $g$ we can assume without loss of generality 
that we consider a function $f'\in \C[a,b,c;u_1,\dots,u_m;x,y]_{abc}$.  Introduce a new variable $z=ay-bc$, 
then $g=abx^2+z^2$, while $f'=\sum_{k,l\geq 0} h_{kl} x^k z^l$, where $h_{kl}\in \C[a,b,c,u_1,\dots,u_m]_{abc}$.
Since the level curves of $\left.(f')^0=f'\right\vert_{a=a^0,b=b^0,c=c^0,u_i=u_i^0}$ and 
$\left.g^0=g\right\vert_{a=a^0,b=b^0,c=c^0,u_i=u_i^0}$ coincide, the Jacobian 
$$
\begin{vmatrix}
{\partial f'}/{\partial x} & {\partial g}/{\partial x}  \\
{\partial  f'}/{\partial z} & {\partial g}/{\partial z} 
\end{vmatrix}
$$
vanishes in an open subset of $\C^{m+5}$. 
Consequently,
$$
z\frac{\partial f'}{\partial x}-abx\frac{\partial f'}{\partial z}=0,
$$
%Then, $\sum_{k,l} z k h_{k,l} x^{k-1} z^l-abx l h_{k,l} x^k z^{l-1}=0$.
and hence, 
$$
\sum_{k,l} \left( (k+1) h_{k+1,l-1} - (l+1)ab h_{k-1,l+1}\right)x^k z^l=0.
$$
Equivalently, for all $k,l$ the equality 
$$
(k+1) h_{k+1,l-1} - (l+1)ab h_{k-1,l+1}=0
$$ 
holds. Solving this recurrence relation in $k$ and 
taking into account that $h_{kl}=0$ for $k<0$ or $l<0$, we see that
$h_{kl}=0$ if at least one of $k$ and $l$ is odd,
while for even $k=2r$, $l=2n-2r$,
$$ 
h_{2r,2n-2r}=(ab)^r{\binom n  r} h_{0,2n}.
$$

Therefore, $f'=\sum_{k=0}^N h_{0,2k} g^k$, and hence  $f=\sum_{k=-N_1}^{N- N_1} h_{0,2k} g^k$,
where $h_{0,2k}\in \C[a,b,c,u_1,\dots,u_m]_{abc}$. 
\qed

\begin{lemma}\label{lem:fR0}
 Let $f\in \wM$, $\ord f=0$, and   
\begin{equation}\label{eq:logcan}
\{x_{12},f\} =\omega x_{12}f.
\end{equation} 
Then $[f]=F(p_0)=\sum_{i=-n_1}^{n_2} a_i p_0^i$, where  
$a_i\in M$.
%\C[x_{11},x_{13},x_{21},S_1,S_2,S_3]_{(x_{11}x_{13}x_{21})}$.
\end{lemma}

{\bf Proof:} By Lemma \ref{lem:decomp}, $\{x_{12},f\} =\{x_{12},[f]\}+x_{12} \{x_{12},\tilde f\}$,
therefore $[\{x_{12},f\}] =[\{x_{12},[f]\}]$, and hence  
\eqref{eq:logcan} yields  $[\{x_{12},[f]\}]=0$. On the other hand,
$[\{x_{12},[f]\}]=[(D_1+D_2)([f])]=D_2([f])$, 
hence $D_2([f])=0$.
Solving the differential equation 
$$
(2x_{13}x_{21}-2x_{11}x_{23})\frac{\partial [f]}{\partial x_{22}}+
(2x_{13} x_{22})\frac{\partial [f]}{\partial x_{23}}=0
$$
by the method of characteristics we obtain a parametric equation 
\begin{equation}\label{eqn:parameter}
\begin{cases}
x_{22}(t)=\sqrt{\frac{\rho}{x_{11}x_{13}}}\cdot\cos(2\sqrt{x_{11}x_{13}}\cdot t),\\
x_{23}(t)=\frac{1}{x_{11}}(\sqrt{\rho}\cdot\sin(2\sqrt{x_{11}x_{13}}\cdot t)+x_{13}x_{21}).
\end{cases}
\end{equation}
 Therefore, characteristics are level curves $E_{\rho}$ described by the equation $p_0=\rho$. 
 Consequently, Lemma~\ref{lem:level} applies with $a=x_{11}$, $b=x_{13}$, $c=x_{21}$, $x=x_{22}$, $y=x_{23}$, 
 $m=3$, $u_i=s_i$ for $i=1,2,3$, and $g=p_0$. According to this lemma,  
 $[f]$ is a Laurent polynomial in $p_0$, $[f]=F(p_0)=\sum_{k=-n_1}^{n_2} a_k p_0^k$ with 
 $a_k\in M$,  which accomplishes the proof.
\qed

\begin{lemma}\label{lem:fR} Let $f$ satisfy all conditions of Lemma \ref{lem:fR0}, 
so that
$[f]=F(p_0)=\sum_k a_k p_0^k$. Then $\{x_{12},F(p)\}=\omega x_{12} F(p)$.
\end{lemma}

{\bf Proof:} By the proof of Lemma~\ref{lem:fR0},
$\{x_{12},[f]\}=D_1([f])$. Note that $D_2={d}/{dt}$, where $t$ is the parameter 
of the level sets $E_{\rho}= \{p_0=\rho\}$ in \eqref{eqn:parameter}. In particular, 
fixing real values of all variables except for $x_{22}$ and $x_{23}$
and assuming $x_{11}x_{13}>0$ and $\rho>0$, the level set  $E_{\rho}$ becomes an ellipse in the 
real plane $x_{22},x_{23}$. Condition $\{x_{12},f\}=\omega x_{12} f$ implies
$D_1([f])+D_2(x_{12}\tilde f)=\omega x_{12} [f]$, or, equivalently, 
$\frac{d}{dt} x_{12}\tilde f = \omega x_{12} [f] -D_1([f])$.
Recall that $x_{12}\tilde f$ is a rational  function in $x_{22}, x_{23}$, which means that 
its integral along any level set $E_{\rho}$ vanishes:  
\begin{equation}\label{zerint}
\int_{E_{\rho}} (\omega x_{12} [f] -D_1([f])) dt=0.
\end{equation}

Assuming $[f]=F(p_0)=\sum_{k,i}  a_{ki} p_0^k$ where each $a_{ki}$ is a Laurent monomial in 
$M$, we get 
$$
D_1([f])=x_{12}\left(\sum_{k,i}  \alpha_{ki} a_{ki} p_0^k+
\sum_{k,i} k  a_{ki} p_0^{k-1} D_1(p_0)
\right)
$$
for some rational constants $\alpha_{ki}$ determined by the equation
$\{x_{12},a_{ki}\}=\alpha_{ki} x_{12} a_{ki}$; by \eqref{coorbra}, such a constant exists 
for any monomial in $M$.

Let us find $D_1(p_0)$:
\begin{equation*}
\begin{aligned}
\frac 3{x_{12}}D_1(p_0)&={8}x_{11}x_{13}x_{22}^2+
{8}\Delta
\left(-\frac{1}{2}x_{11}x_{23}-x_{13}x_{21}\right)\\
&={8}p_0-12\Delta^2-12\Delta x_{13}x_{21}
\end{aligned}
\end{equation*}
with $\Delta=x_{11}x_{23}-x_{13}x_{21}$.

Note that 
$$
\int_{E_{\rho}}p_0^kdt=\pi\rho^k/{\sqrt{x_{11}x_{13}}}
$$
since $p_0$ equals $\rho$ on $E_{\rho}$, 
$$
\int_{E_{\rho}}4(x_{11}x_{23}-x_{13}x_{21})x_{13}x_{21} dt=0
$$ 
by the symmetry of the level set $E_{\rho}$, and  
$$
\int_{E_{\rho}}4(x_{11}x_{23}-x_{13}x_{21})^2 dt=2\pi \rho/\sqrt{x_{11}x_{13}}.
$$
Therefore, \eqref{zerint} yields
\[
\begin{aligned}
&\frac 1{x_{12}}\int_{E_{\rho}} (\omega x_{12} [f] -D_1([f])) dt\\
&=\frac\pi{\sqrt{x_{11}x_{13}}}\left(\omega \sum_{k,i} a_{ki} \rho^k  - 
\sum_{k,i} \alpha_{ki} a_{ki} \rho^k
-\sum_{k,i} \frac{8}{3} k  a_{ki} \rho^k \right)\\
&+\sum_{k,i} \frac{2\pi \rho}{\sqrt{x_{11}x_{13}}} k  a_{ki} \rho^{k-1}=0.
\end{aligned}
\]
Equivalently, 
$$
\frac{\pi}{\sqrt{x_{11}x_{13}}}a_{kI} \rho^k \left(\omega -\alpha_{ki} -\frac{2}{3}k\right)=0,
$$
or 
\begin{equation}\label{eq:c}
\omega=\alpha_{ki}+\frac{2}{3}k
\end{equation} 
for all $k$ and $i$.

%In particular, we note that all $\alpha_{kI}$ are equal for the same $k$ and different $I$,
%and we can introduce $\alpha_k=\alpha_{kI}$ for any $I$ such that all 

Let us compute the Poisson bracket of $x_{12}$ and  $F(p)=\sum_k \sum_i a_{ki} p^k$. 
Note that $\{x_{12},p\}=\frac{2}{3}x_{12} p$, hence
$$
\left\{x_{12},\sum_{k,i} a_{ki} p^k\right\} = 
x_{12}\sum_{k,i} \left(\alpha_{ki} +\frac{2}{3} k\right) a_{ki} p^k.
$$ 
Taking into account \eqref{eq:c}, we conclude that 
$\{x_{12},F(p)\}=\omega x_{12} F(p)$.
 \qed 

\begin{corollary}\label{cor:approx} 
 Let $f\in \wM$  and $\{x_{12},f\}=\omega x_{12} f$.
Then for any positive integer $d>0$ there exists a Laurent polynomial $F_d(\xi)=\sum_k a_{dk} \xi^k$ such that 
$a_{dk}\in M$, $\{x_{12},F_d(p)\}=\omega x_{12}F_d(p)$ and $\ord (f-F_d(p))\geq d$.
\end{corollary}

{\bf Proof:} The proof goes by induction on $d$. For $d=1$, define $f_1=f/x_{12}^{\ord f}$. Then $f_1$ satisfies 
all conditions of Lemma \ref{lem:fR0}, and hence $[f_1]=F_1(p_0)$. Therefore, by Lemma \ref{lem:fR}, 
$F_1(p)$ possesses all the desired properties.

Assume that we have built $F_{d-1}$. Define $f_d=(f-F_{d-1}(p))/x_{12}^{\ord(f-F_{d-1}(p))}$. Clearly, $\ord f_d=0$,
and by the inductive hypothesis, $f_d$ satisfies \eqref{eq:logcan}. Therefore, $[f_d]=F(p_0)$ and 
$\ord (f_d-F(p))\geq 1$. Consequently, it suffices to put $F_d(p)=F_{d-1}(p)+x_{12}^{\ord(f-F_{d-1}(p))}F(p)$.  
\qed

\begin{lemma}\label{lem:twofunctions} Let 
$f,g\in \wM$, and assume that $\{x_{12},f\}=\omega_1 x_{12} f$ and
$\{x_{12},g\}=\omega_2 x_{12} g$.
Then
$$
J=\begin{vmatrix}
{\partial f}/{\partial x_{22}} & {\partial g}/{\partial x_{22}} \\
{\partial f}/{\partial x_{23}} & {\partial g}/{\partial x_{23}} 
\end{vmatrix}=0.
$$
\end{lemma}

{\bf Proof:}
Assume that $J\ne 0$, then $d=\ord J \ge 0$ is finite. 
%We prove the statement by contradiction showing that $\ord J>d$. Indeed, 
By Corollary~\ref{cor:approx}, there exist Laurent polynomials
$F_d$ and $G_d$ such that $f=F_d(p)+\widetilde{F}_{d}$, $\ord \widetilde{F}_{d}>d$ and 
$g=G_d(p)+\widetilde{G}_{d}$, $\ord \widetilde{G}_{d}>d$. Therefore,
\[
\begin{aligned}
J&=\begin{vmatrix}
\dfrac{\partial F_d(p)}{\partial x_{22}}+\dfrac{\partial \widetilde{F}_{d}}{\partial x_{22}}& 
\dfrac{\partial G_d(p)}{\partial x_{22}} +\dfrac{\partial \widetilde{G}_{d}}{\partial x_{22}}\\
\dfrac{\partial F_d(p)}{\partial x_{23}}+\dfrac{\partial \widetilde{F}_{d}}{\partial x_{23}} & 
\dfrac{\partial G_d(p)}{\partial x_{23}} +\dfrac{\partial \widetilde{G}_{d}}{\partial x_{23}}
\end{vmatrix}\\
&=\begin{vmatrix}
\dfrac{\partial F_d(p)}{\partial x_{22}}& \dfrac{\partial G_d(p)}{\partial x_{22}} \\
\dfrac{\partial F_d(p)}{\partial x_{23}}& \dfrac{\partial G_d(p)}{\partial x_{23}} 
\end{vmatrix}+\widetilde{J}
\end{aligned}
\]
with $\ord\widetilde{J}>d$. It remains to notice that
{\renewcommand{\arraystretch}{2.5}
\[
 \begin{vmatrix}
\dfrac{\partial F_d(p)}{\partial x_{22}}& \dfrac{\partial G_d(p)}{\partial x_{22}} \\
\dfrac{\partial F_d(p)}{\partial x_{23}}& \dfrac{\partial G_d(p)}{\partial x_{23}} 
\end{vmatrix}=
 \begin{vmatrix}
\dfrac{dF_d}{dp}\dfrac{\partial p}{\partial x_{22}}& \dfrac{dG_d}{dp}\dfrac{\partial p}{\partial x_{22}} \\
\dfrac{dF_d}{dp}\dfrac{\partial p}{\partial x_{23}}& \dfrac{dG_d}{dp}\dfrac{\partial p}{\partial x_{23}} 
\end{vmatrix}=0,
\]}
and hence $\ord J=\ord\widetilde{J}>d$, 
a contradiction.
\qed

\begin{corollary}\label{cor:algdep} 
The function $x_{12}$ is not a cluster variable in the cluster algebra $\A_{CG}(3)$.
\end{corollary}

{\bf Proof:}
Assume that there exists a cluster containing $x_{12}$. 
This cluster contains also three stable variables $s_1, s_2, s_3$ and five cluster variables  
$v_1, v_2, v_3, v_4, v_5\in\C[\Mat_3]$ such that $\{x_{12},v_j\}=\omega_j x_{12} v_j$. 
%and the cluster variables form algebraically inmdependent $9$-tuple. 
Choose $x_{11}, x_{12}, x_{13}, x_{21}, x_{22}, x_{23}, s_1, s_2, s_3$ as local coordinates 
in an open subset of $\C^9$ and  
write all cluster variables as elements of  $\wM$.
Since the cluster variables belonging to the same cluster are algebraically independent,  the rank 
of the $2\times 5$ matrix
$$
\Jac=\begin{pmatrix}
\dfrac{\partial v_1}{\partial x_{22}} & \dfrac{\partial v_2}{\partial x_{22}} & 
\dfrac{\partial v_3}{\partial x_{22}} & \dfrac{\partial v_4}{\partial x_{22}} & 
\dfrac{\partial v_5}{\partial x_{22}} \\
\dfrac{\partial v_1}{\partial x_{23}} & \dfrac{\partial v_2}{\partial x_{23}} & 
\dfrac{\partial v_3}{\partial x_{23}} & \dfrac{\partial v_4}{\partial x_{23}} & 
\dfrac{\partial v_5}{\partial x_{23}}
\end{pmatrix}
$$
equals 2 for generic values of parameters  $x_{11}, x_{12}, x_{13}, x_{21}, s_1, \allowbreak s_2, s_3$. 
However, by Lemma~\ref{lem:twofunctions}, any $2\times 2$ minor of $\Jac$ vanishes identically, a contradiction.
\qed

 \begin{corollary}\label{cor:notinclusteralgebra} 
The function $x_{12}$ does not belong to the cluster algebra $\A_{CG}(3)$.
\end{corollary}

{\bf Proof:}
Note that cluster algebra $\A_{CG}(3)$ has three independent compatible toric actions:
$T_1(x_{ij})=t^{i-2}\cdot x_{ij}$, $T_2(x_{ij})=t^{j-2}\cdot x_{ij}$, and $T_3(x_{ij})=t\cdot x_{ij}$. 
Any cluster variable of $\A_{CG}(3)$ is a polynomial in $x_{ij}$ homogeneous with respect to all three weights $w_1$, $w_2$, $w_3$ determined by the toric actions $T_1$, $T_2$, $T_3$.
Therefore, we can associate  with every cluster coordinate $z_C$ in a cluster $C$ the weight vector  $w(z_C)=(w_1(z_C),w_2(z_C))$ containing only the first two weights. 
Note that the weight vectors for all nine matrix entries $x_{ij}$, $1\le i,j\le 3$, are all distinct, implying that no cluster variable is a linear combination of matrix entries $x_{ij}$ with 
more than one nontrivial coefficient.

Taking that into account and recalling that the cluster algebra contains only regular function on $\Mat_3$,
 we see that if any linear polynomial $\sum c_{ij} x_{ij}$ belongs to $\A_{CG}(3)$  then
$c_{ij}\ne 0$ if and only if $x_{ij}$ itself is a cluster variable in $\A_{CG}(3)$.    
Since $x_{12}$ is not a cluster variable by Corollary~\ref{cor:algdep}, we conclude that
$x_{12}\notin \A_{CG}(3)$.
\qed

\begin{remark}
The same arguments imply that $x_{12}$ does not belong to the Cremmer--Gervais cluster algebra 
$\A'_{CG}(3)$ associated with $SL_3$. Note first that again $x_{12}$ is not a cluster variable in $\A'_{CG}(3)$.
 Indeed, if there is a cluster $C$ in $\A'_{CG}(3)$ containing $x_{12}$, then
by adding $\det X$ to $C$ we will obtain a cluster in $\A_{CG}(3)$ containing $x_{12}$, in a contradiction 
with Corollary~\ref{cor:algdep}. The proof of Corollary~\ref{cor:notinclusteralgebra} can be used literally to conclude that  $x_{12}\notin \A'_{CG}(3)$.
\end{remark}

{\section{Connection to total positivity}} 
 
A comprehensive analysis of total positivity in reductive groups that was performed in \cite{FZ1}
gave a major impetus to the development of the theory of cluster algebras. As it was explained 
in \cite[Remark 2.16]{CAIII}, the set of totally positive elements in a simple Lie group coincides with the set of elements in an open double Bruhat cell that form a positive locus with respect to the standard cluster structure. This locus is formed by elements on which  functions that form a single extended cluster
(and therefore all cluster variables) take positive values. 

In the case of $GL_n$, the set $TP(n)$ of totally positive elements is formed by matrices with all minors positive. It is a well-studied class of matrices that has many important applications. It is natural to ask what can be said about the positive locus in $GL_n$ with respect to the Cremmer--Gervais cluster structure introduced above. We denote this locus by $TP_{CG}(n)$. A careful analysis of cluster transformations that form the sequence $\mathcal T$  featured in Theorem \ref{transform2} allows us to make the first step towards understanding $TP_{CG}(n)$.

\begin{theorem}
\label{TP}
$$
TP_{CG}(n) \subsetneq TP(n).
$$
\end{theorem}

{\bf Proof.} It is shown in \cite{GSV6} that in the process of applying transformations forming $\mathcal T$ to the initial cluster described in Theorem \ref{logcan}, one recovers the following two families
of minors as cluster variables in $\A_{CG}$: 
$$
F_1=\{ \det X^{[i,n]}_{[j,n+j-i]}\: i\in [n], j\in [i-1] \}
$$
%\footnote{ \cite{GSV6}, Lemma 7.9, line 2 : $l=i-j, p=i-1, q=i-j+1$} 
and 
$$
F_2=\{ \det X^{[1,n]\setminus [j+1,j+n-i]}_{[1,i]}\: i\in [n], j\in [i] \}.
$$
%\footnote{ \cite{GSV6}, Lemma 7.12, line 3 : $l=j, p=n+j-i, q=n-i$}. 
By Theorem \ref{transform2}, the family
\begin{align*}
F^{w_0}_1&=\{ \det (W_0 X W_0)^{[i,n]}_{[j,n+j-i]}\\
&= \det X ^{[1,n-i+1]}_{[i-j+1,n-j+1]}\: i\in [n], j\in [i-1] \} 
\end{align*}
also consists of cluster variables in $\A_{CG}$. But according to \cite[Thm.~4.13]{FZ1}, $F^{w_0}_1 \cup F_2$ is one of the {\em test families\/} for total positivity
(it corresponds to the reduced word $\overline {n-1} \ldots \bar {1}\overline {n-1} \ldots \bar {2}\ldots \overline {n-1} 1\ldots n-1 1 \dots n-2 \ldots 1$), 
i.e. positivity of all minors in this family guarantees total positivity of $X$. This shows that $TP_{CG}(n) \subset TP(n)$.

To complete the proof, we need to construct a totally positive matrix $X$ such that at least one of the functions forming the initial cluster for $\A_{CG}$ is negative when evaluated at $X$. Let $S=(\delta_{i,j-1})_{i,j=1}^{n}$ be the $n\times n$ shift matrix. For $t\in \C$, define $E(t)= \one + t S$. It is easy to see that if  $t > 0$ then
matrices $E(t)$ and $E(-t)^{-1}$ are totally nonnegative. Let $X_0= \left ( E(-1)^T E(-1)\right )^{-1}$ and
$X(t) = X_0 \left (E(t) E(t)^T\right )^{n-1}$. The Cauchy--Binet identity implies that $X_0$ is totally nonnegative and  $X(t)= X_0 (\one + O(t))$ is totally positive if $t > 0$. However, the value of the  cluster variable
$\psi_2(X) = x_{n-2,n} x_{n 1} - x_{n-1,n} x_{n-1, 1}$ at $X_0$ is $-1$ and so, for $t$ positive and sufficiently small, $X(t)\in TP(n)\setminus TP_{CG} (n)$.
\qed

\section{Further directions}

We would like to outline briefly further directions of research.

1. Our main goal is to prove Conjecture~\ref{ulti}. In addition to the results reported above,
the conjecture has been verified for $SL_5$, see~\cite{Idan}. In all the cases studied so far, one can produce
a quiver that is remarkably similar to the one shown in Fig.~1 above. Namely, it consists of $n^2$ vertices placed on a grid with diagonals. Besides,  vertices of the upper and the lower rows, as well as vertices
of the leftmost and the rightmost columns are connected with additional edges, which are determined by the isometry
$\gamma$. 

2. We conjecture that the result of Theorem~\ref{gap} can be extended to other exotic cluster structures as well. 

3. Theorems~\ref{logcan} and~\ref{regular} follow from similar results valid in the double. This hints at a possibility of endowing the double with a cluster structure associated with the Cremmer--Gervais bracket.

\section*{Acknowledgments}
 M.~G.~was supported in part by NSF Grants DMS \#0801204 and  DMS \# 1101462. 
M.~S.~was supported in part by NSF Grants DMS \#0800671 and DMS \# 1101369.  
A.~V.~was supported in part by ISF Grant \#162/12. This paper was partially written during
our joint stays at MFO Oberwolfach (Research in Pairs program, August 2010), at the Hausdorff Institute
(Research in Groups program, June-August, 2011) and at the MSRI (Cluster Algebras program, August-December, 2012). We are grateful to these institutions for
warm hospitality and excellent working conditions. We are also grateful to our home institution for support in arranging visits by collaborators.
Special thanks are due to Bernhard Keller, whose quiver mutation applet proved to be an indispensable tool in our work on this project, to Thomas Lam for attracting our attention to the total positivity question, and to Sergei Fomin for illuminating discussions.


\begin{thebibliography}{O}

\bibitem{BD} A.~Belavin and V.~Drinfeld,
\textit{Solutions of the classical Yang-Baxter equation for simple Lie algebras}.
Funktsional. Anal. i Prilozhen. {\bf16} (1982), 1-�29.
%\textit{Triangle equations and simple Lie algebras}.
%Soviet Sci. Rev. Sect. C Math. Phys. Rev. {\bf 4} (1984), 93--165. 

\bibitem{CAIII}  A.~Berenstein, S.~Fomin, and A.~Zelevinsky,
\textit{Cluster algebras. III. Upper bounds and double Bruhat cells}. 
Duke Math. J. \textbf{126} (2005), 1--52.

%\bibitem{CP} V. Chari, A. Pressley, 
%\textit{A guide to quantum groups}. 
%Cambridge University Press, 1994.

%\bibitem[BGY]{BGY} K. A.~Brown, K.  R.~Goodearl and M.~Yakimov, \textit{Poisson structures on affine spaces and flag varieties. I. Matrix affine Poisson space}, Adv. Math. {\bf 206} (2006), no. 2, 567--629.

%\bibitem[D]{D} V. G.~ Drinfeld, \textit{On Poisson homogeneous spaces of Poisson--Lie groups}, Theoret. and Math. Phys. {\bf 95} (1993), no. 2, 524--525.

\bibitem{Idan} I.~Eisner,
\textit{Exotic cluster structures on $SL_5$}.
in preparation.

\bibitem{FZ1} S.~Fomin and A.~Zelevinsky,  
\textit{Double Bruhat cells and total positivity}. 
J. Amer. Math. Soc. \textbf{12} (1999), 335-–380.

\bibitem{FZ2}  S.~Fomin and A.~Zelevinsky, 
\textit{Cluster algebras.I. Foundations}. 
J. Amer. Math. Soc. \textbf{15} (2002), 497--529.

\bibitem{FZ3} S.~Fomin and A.~Zelevinsky, 
\textit{The Laurent phenomenon}.
Adv. in Appl. Math. {\bf 28} (2002), 119--144.

\bibitem{GSV1}  M.~Gekhtman, M.~Shapiro, and A.~Vainshtein,
\textit{Cluster algebras and Poisson geometry}.  
Mosc. Math. J. \textbf{3} (2003), 899--934.

\bibitem{GSV5}  M.~Gekhtman, M.~Shapiro, and A.~Vainshtein, 
\textit{Generalized B\"acklund-Darboux transformations of Coxeter-Toda flows from a cluster algebra perspective}
Acta Math. \textbf{206} (2011), 245--310.

\bibitem{GSVb}  M.~Gekhtman, M.~Shapiro, and A.~Vainshtein,
\textit{Cluster algebras and Poisson geometry}.
Mathematical Surveys and Monographs, 167. American Mathematical Society, Providence, RI, 2010.


\bibitem{GSVMMJ}  M.~Gekhtman, M.~Shapiro, and A.~Vainshtein,
\textit{Cluster structures on simple complex Lie groups and Belavin--Drinfeld classification}, Mosc. Math. J. \textbf{12} (2012), 293--312.

\bibitem{GSV6} M.~Gekhtman, M.~Shapiro, and A.~Vainshtein,
\textit{Exotic cluster structures in $SL_n$: the Cremmer--Gervais case}, arXiv:1307.1020 

%\bibitem{GeGi}  M.~Gerstenhaber and  A.~Giaquinto, 
%\textit{Boundary solutions of the classical Yang-Baxter equation}, 
%Lett. Math. Phys. {\bf 40} (1997),  337--353. 

%\bibitem{KoSo}  L.~Korogodski and  Y.~Soibelman, 
%\textit{Algebras of functions on quantum groups. Part I}. 
%Mathematical Surveys and Monographs, 56. American Mathematical Society, Providence, RI, 1998.

\bibitem{r-sts}  A.~Reyman and M.~Semenov-Tian-Shansky
\textit{Group-theoretical methods in the theory of
finite-dimensional integrable systems} Encyclopaedia of
Mathematical Sciences, vol.16, Springer--Verlag, Berlin, 1994 pp.116--225.

\bibitem{Ya} M.~Yakimov, 
\textit{Symplectic leaves of complex reductive Poisson-Lie groups}.
Duke Math. J. {\bf 112} (2002), 453--509.
\end{thebibliography}
\end{document}